%% file: main.tex
\documentclass{IEEEtran4PSCC}
\ifCLASSINFOpdf
   \usepackage[pdftex]{graphicx}
\else
   \usepackage[dvips]{graphicx}
\fi
\usepackage[cmex10]{amsmath}
\usepackage{amssymb}
\newcommand{\II}{\mathrm{II}}

\newcommand{\T}{\top}
\newcommand{\gausstail}[1]{\Phi\!\left(-#1\right)}

\makeatletter
\let\old@ps@headings\ps@headings
\let\old@ps@IEEEtitlepagestyle\ps@IEEEtitlepagestyle
\def\psccfooter#1{%
    \def\ps@headings{%
        \old@ps@headings%
        \def\@oddfoot{\strut\hfill#1\hfill\strut}%
        \def\@evenfoot{\strut\hfill#1\hfill\strut}%
    }%
    \def\ps@IEEEtitlepagestyle{%
        \old@ps@IEEEtitlepagestyle%
        \def\@oddfoot{\strut\hfill#1\hfill\strut}%
        \def\@evenfoot{\strut\hfill#1\hfill\strut}%
    }%
    \ps@headings%
}
\makeatother

\psccfooter{%
        \parbox{\textwidth}{\hrulefill \\ \small{24th Power Systems Computation Conference} \hfill \begin{minipage}{0.2\textwidth}\centering \vspace*{4pt} \includegraphics[scale=0.06]{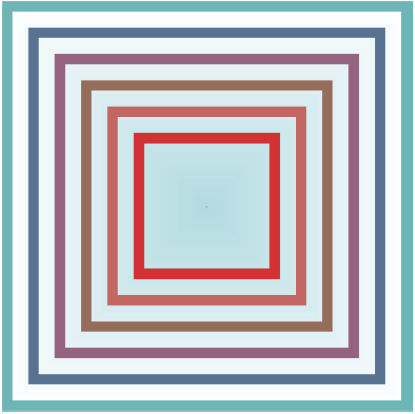}\\\small{PSCC 2026} \end{minipage} \hfill \small{Limassol, Cyprus --- June 8-12, 2026}}%
}

\usepackage{xcolor}

\allowdisplaybreaks

\begin{document}
\title{Computing Rare Probabilities of Voltage Collapse}

\author{\IEEEauthorblockN{Tongtong Jin\IEEEauthorrefmark{1},
Anirudh Subramanyam\IEEEauthorrefmark{2},
D. Adrian Maldonado \IEEEauthorrefmark{1} \IEEEauthorrefmark{3}}
\IEEEauthorblockA{\IEEEauthorrefmark{1} Computational and Applied Mathematics \\
The University of Chicago,
5801 S Ellis Ave, Chicago, IL 60637}
\IEEEauthorblockA{\IEEEauthorrefmark{2} Department of Industrial and Manufacturing Engineering\\
The Pennsylvania State University,
201 Old Main, University Park, PA 16802}
\IEEEauthorblockA{\IEEEauthorrefmark{3} Mathematics and Computer Science Division\\
Argonne National Laboratory,
9700 S Cass Ave, Lemont, IL 60439}
}

%% To specify the authors when (number of affiliations > 2)
% \author{\IEEEauthorblockN{Author n.1\IEEEauthorrefmark{1},
% Author n.2\IEEEauthorrefmark{2},
% Author n.3\IEEEauthorrefmark{3}, 
% Author n.4\IEEEauthorrefmark{3} and
% Author n.5\IEEEauthorrefmark{4}}
% \IEEEauthorblockA{\IEEEauthorrefmark{1} Department Name of Organization A\\
% Name of the organization A,
% Address A\\ Emails if wanted}
% \IEEEauthorblockA{\IEEEauthorrefmark{2} Department Name of Organization B\\
% Name of the organization B,
% Address B\\ Emails if wanted}
% \IEEEauthorblockA{\IEEEauthorrefmark{3} Department Name of Organization C\\
% Name of the organization C,
% Address C\\ Emails if wanted}
% \IEEEauthorblockA{\IEEEauthorrefmark{4}Department Name of Organization D\\
% Name of the organization D,
% Address D\\ Emails if wanted}
% }

\maketitle

\begin{abstract}
This paper introduces a framework based on Large Deviation Theory (LDT) to accurately and efficiently compute the rare probabilities of voltage collapse. We formulate the problem as finding the most probable failure point (the instanton) on the stability boundary and derive both first-order and second-order approximations for the collapse probability. The second-order method incorporates the local curvature of the stability boundary, yielding higher accuracy. This LDT framework generalizes methods based on Mahalanobis distance and is extensible to non-Gaussian uncertainties. We validate our approach on test systems, demonstrating that the LDT estimates converge to Monte Carlo results in the rare-event regime where direct sampling becomes computationally prohibitive.
\end{abstract}


\thanksto{This material is based upon work supported by the National Science Foundation under Grant No. DMS-2229408 and DMS-2229409.}

\section{Introduction}

Maintaining voltage stability is critical for preventing high-impact, rare events like cascading blackouts. Foundational research in this area established a deterministic understanding of stability through geometric margins. The influential work of Dobson and Lu \cite{dobson1993}, for instance, provided a rigorous method for computing the closest saddle-node bifurcation (SNB), defining an essential, worst-case distance to instability.

Building on this, Zhang et al. \cite{ZHANG2004} introduced a probabilistic approach to quantify the Transmission Reliability Margin (TRM). Their method uses first-order sensitivities of the transfer capability with respect to various uncertain parameters. By assuming a large number of independent uncertainties, they invoke the Central Limit Theorem to approximate the total uncertainty in the stability margin as a normal distribution. This represents a key step towards probabilistic assessment but is fundamentally a first-order, linear approximation that may lose accuracy when faced with large uncertainties or highly nonlinear stability boundaries.

More recently, the work of Mittelstaedt et al. \cite{Mittelstaedt2017} provided a crucial link between the geometric and probabilistic viewpoints. They proposed an iterative method to find the most probable bifurcation point by minimizing the Mahalanobis distance from the mean operating point to the stability boundary. For Gaussian uncertainties, minimizing the Mahalanobis distance is mathematically equivalent to finding the point that maximizes the probability density. Their approach identifies the most likely failure point but does not provide an estimate of the total collapse probability.

A key challenge in applying probabilistic methods to voltage collapse is the rare-event nature of such occurrences. Accurately computing these low probabilities can be computationally demanding or inaccurate with standard methods. Techniques such as large deviation theory (LDT) attempt to solve this issue. 
The core intuition of LDT is that the probability of a rare failure is overwhelmingly dominated by the single most probable failure scenario, known as the "instanton". Rather than sampling randomly, LDT seeks to find this specific point directly through optimization. 
In addition to estimating the probability of voltage collapse, the resulting instanton identifies the most likely collapse-inducing loading pattern, thereby providing operational insight into the dominant direction of stress and potential preventive actions.

In power systems analysis, Chertkov et al. \cite{Chertkov2011} have applied this instanton-based approach to identify the most probable load patterns that cause static line overloads and require load shedding within a simplified DC power flow model.
Nesti et al. \cite{nesti2019temperature} develop analytical LDT approximations for rare line overheating probabilities under Ornstein-Uhlenbeck injection models.
Related works use other LDT-based line tripping probabilities to develop cascading failure simulators \cite{roth2021kinetic}, to mitigate cascade risk during generation dispatch \cite{subramanyam2022failure}, and for network vulnerability analysis \cite{lam2023network}.
More recently, Tapia et al. \cite{tapia2024electricity} applied the LDT-based probabilistic optimization framework of Tong et al. \cite{tong2022optimization} to introduce a new extreme reserve product in electricity market clearing.

Our work fills an important gap in the literature on voltage collapse by providing a method to accurately and efficiently compute these low-probability events. In particular, our contributions are

\begin{itemize}
    \item We formulate voltage-collapse risk as a rare-event probability problem in the framework of large deviation theory (LDT). This yields an instanton-based characterization of the most probable collapse scenario together with explicit first-order approximations of the collapse probability. In this way, earlier Mahalanobis-distance-based approaches such as \cite{Mittelstaedt2017} are recovered as the Gaussian special case, while the same framework extends naturally to non-Gaussian uncertainty models such as Gaussian mixtures.
    \item We derive a second-order approximation of the collapse probability that incorporates the local geometry of the bifurcation boundary at the instanton. The resulting curvature-dependent correction can be computed efficiently in closed form and improves probability estimation when the stability boundary is not locally well approximated by its tangent hyperplane.

\end{itemize}

\section{Methodology}

Voltage collapse is modeled as the disappearance of a stable operating equilibrium in response to changing system parameters. While the collapse itself is a dynamic event, the boundary of stability can be analyzed using the system's static power flow equations \cite{dobson2011irrelevance}
\begin{equation}\label{eq:pfe}
    f(x, \lambda) = 0
\end{equation}
where $x \in \mathbb{R}^n$ is the vector of state variables (e.g., voltage magnitudes and angles) and $\lambda \in \mathbb{R}^m$ is a vector of uncertain parameters, such as deviations from some mean load and generation level. At a stable operating point, the Jacobian matrix, $f_x(x, \lambda) = \frac{\partial f}{\partial x}$, must be non-singular.

Voltage collapse corresponds to the disappearance of a stable operating equilibrium, an event mathematically described as a saddle-node bifurcation \cite{dobson2011irrelevance}. The primary condition for this is that the Jacobian matrix $f_x$ becomes singular. In addition, transversality and non-degeneracy conditions must be satisfied to distinguish a generic fold bifurcation from other, more complex types. As these conditions are typically met in power system models, the set of points defined by a singular Jacobian is reliably taken to be the stability boundary for voltage collapse.

The set of all such points forms the stability boundary, or bifurcation surface, in the parameter space:

\begin{align}
\mathcal{B} 
&= \{ \lambda \in \mathbb{R}^m : \text{there exists } x \text{ such that } f(x,\lambda)=0 \nonumber \\
&\qquad \text{and } \det(f_x(x,\lambda))=0 \}.
\end{align}
This surface $\mathcal{B}$ is typically a smooth $(m-1)$-dimensional manifold that divides the parameter space into a feasible (stable) region $\mathcal{B}^{+}$ and an infeasible region $\mathcal{B}^-$ where no stable equilibrium exists \cite{dobson2011irrelevance}. As shown in \cite{appliedbifurcation}, the saddle-node (fold) bifurcation is a manifold of codimension 1 in the parameter space.
Our objective is to calculate the probability that a random vector $\lambda$ lies inside this infeasible region.

\subsection{Large Deviation Theory}

Large Deviation Theory (LDT) provides an asymptotic estimate for the probability of rare events \cite{dembo2009large,Bucklew_2013B}. 
Intuitively, the theory states that for a ``rare'' set $\mathcal{A}$, the probability $P(\lambda \in \mathcal{A})$ can be approximated as follows:
\begin{equation}\label{eq:ldp_approximation}
P(\lambda \in \mathcal{A}) \approx e^{-I(\lambda^*)},
\end{equation}
where $I(\cdot)$ is the so-called rate function and $\lambda^*$ is the point that minimizes this function over $\mathcal{A}$:
\begin{equation}\label{eq:instanton}
I(\lambda^*) = \inf_{\lambda \in \mathcal{A}} I(\lambda).
\end{equation}

Formally, the approximation in~\eqref{eq:ldp_approximation} can be made precise under an appropriate asymptotic regime.
For example, let $\lambda_n := \lambda/\sqrt{n}$ be a scaled-down version of some zero-mean random deviation $\lambda$ from a fixed nominal parameter vector.
If $0 \notin \mathcal{A}$ (e.g., the nominal parameters induce a stable operating point), then observe that, as $n$ increases, the event $\{\lambda_n \in \mathcal{A}\}$ becomes increasingly rare.
In fact, this can be made precise as follows:
\[
\lim_{n \to \infty} -\frac{1}{n} \log P(\lambda_n \in \mathcal{A}) = I(\lambda^{*}),
\]
which directly motivates the approximation~\eqref{eq:ldp_approximation}.
We refer the reader to \cite{burenev2025introduction} for an overview of LDT and its applications.

The point $\lambda^*$ is also called the \textit{instanton}.
Intuitively, it represents the most probable realization of the random parameters that leads to the rare event.
For many distributions (e.g., light-tailed and possessing a density), the rate function is simply the convex conjugate (or Legendre--Fenchel transform) of its cumulant generating function \cite{burenev2025introduction}.

For common distributions, it often has a simple form.
For a multivariate normal $\lambda \sim \mathcal{N}(\mu, \Sigma)$, the rate function coincides with the negative log-density:
\begin{equation}\label{eq:rate_function_mvnormal}
I(\lambda) = \frac{1}{2} (\lambda - \mu)^\top \Sigma^{-1} (\lambda - \mu).
\end{equation}
The minimization in~\eqref{eq:instanton} is therefore equivalent to finding the point on the stability boundary that is closest to the mean $\mu$ in the Mahalanobis sense, which aligns with the intuition in~\cite{Mittelstaedt2017}.
However, the LDT formalism enables us to also consider non-Gaussian distributions for which the rate function can be readily calculated or evaluated, including exponential \cite{burenev2025introduction} and Gaussian mixture distributions \cite{tong2022optimization}; the latter is developed explicitly in Section~\ref{sec:gmm_extension}.

The LDT estimate~\eqref{eq:ldp_approximation} is only asymptotic and hence, rather coarse in practice. Here in the context of power systems and voltage collapse, we consider the collapse region $\mathcal{B}^-$ as $\mathcal{A}$. We shall show in Sections~\ref{sec:1st_order} and~\ref{sec:2nd_order} how information about the first- and second-order geometry of the bifurcation surface $\mathcal{B}$ at the instanton $\lambda^{*}$ can be used to obtain more accurate estimates. The explicit formulas in these sections assume Gaussian uncertainty; the extension to Gaussian mixtures is presented in Section~\ref{sec:gmm_extension}.
We first discuss how to compute the instanton itself.

\subsection{Computing the Instanton}
The instanton $\lambda^*$ is the point on the bifurcation surface $\mathcal{B}$ that minimizes the rate function $I(\lambda)$. This can be formulated as a constrained optimization problem. The first-order optimality (KKT) conditions for this problem lead to the requirement that the gradient of the rate function, $\nabla_\lambda I(\lambda^*)$, must be normal to the bifurcation surface at the instanton.

Similar to the equations derived in \cite{dobsondirect1991} for minimizing $\|\lambda-\lambda_0\|^2$, the conditions for minimizing $I(\lambda)$ can be represented by the following system of equations, which jointly solves for the state $x$, the instanton $\lambda$, the left null vector $w$ of the Jacobian, and a scaling factor $k$:
\begin{subequations} \label{eq:kkt_system}
\begin{align}
f(x, \lambda) &= 0 \,, \\ 
w^\top f_x(x, \lambda) &= 0 \,, \\
k \nabla I(\lambda)^\top - w^\top f_\lambda &= 0 \label{eq:general_optimality} \,, \\
w^\top f_\lambda (w^\top f_\lambda)^\top - 1 &= 0.
\end{align}
\end{subequations}
The solution of this system is denoted as $(x^*, \lambda^*, w^*, k^*)$.

When $\lambda \sim \mathcal{N}(\mu, \Sigma)$ is multivariate normal, the gradient of its rate function~\eqref{eq:rate_function_mvnormal} can be calculated in closed-form:

\begin{equation}\label{eq:rate_gradient}
    \nabla I(\lambda) = \Sigma^{-1} (\lambda - \mu),
\end{equation}
which can then be substituted into condition \eqref{eq:general_optimality}.

When $\lambda$ follows a Gaussian mixture distribution, even though $\nabla I(\lambda)$ does not admit a closed-form expression, it can nevertheless be computed by solving an augmented system of equations involving the cumulant generating function \cite{tong2022optimization} (see Section~\ref{sec:gmm_extension}).
In any case, the full system of nonlinear equations \eqref{eq:kkt_system} can be efficiently solved using Newton-based methods to find the instanton.

\subsection{Normal Vector and Curvature}\label{sec:normal_vector}

To refine the LDT approximation with second-order information, we need geometric information about the bifurcation surface at the instanton: the direction of the outward normal, and the principal curvatures that describe how rapidly the normal vector changes when moving tangentially along the surface. These quantities will appear in the correction term of the second-order LDT approximation (Section~\ref{sec:2nd_order}). We now describe how to compute them.

From \cite{dobson1993computing}, we know that the normal vector to the bifurcation surface at $\lambda$ is:
\begin{equation}
    N=w^\top f_{\lambda} \,.
\end{equation}

The principal curvatures can be computed by taking the eigenvalues of the second fundamental form. The second fundamental form $\II$ is defined by the derivative of the normal vector $N$ with respect to the parameters $\lambda$, denoted by $D_{\lambda}N$. Therefore, to derive the second fundamental form of the bifurcation surface at $\lambda^*$, we need to calculate
\begin{equation} \II = \frac{d}{d\lambda}(w^\top f_{\lambda}) = w^\top f_{xx}(x_{\lambda},x_{\lambda})+w^\top f_{x\lambda}x_{\lambda}+w^\top f_{\lambda\lambda} \,,
\end{equation}
where $f_{xx}(x_{\lambda},x_{\lambda}$) is a bilinear map computed with tensor $f_{xx}$.
In our formulation of the power flow problem, both $f_{x\lambda}$ and $f_{\lambda\lambda}$ are zero \cite{dobson1993computing}. Hence, the second fundamental form can be simplified to:
\begin{equation}\II = w^\top f_{xx}(x_{\lambda},x_{\lambda}) \,.
\end{equation}

To compute $x_{\lambda}$, we then make use of the implicit function theorem as follows:
\begin{subequations}
\begin{align}
f_xx_{\lambda}+f_{\lambda} &= 0 \,, \\
w^\top f_{xx}(x_{\lambda}, v) & = 0 \,,
\end{align}
\end{subequations}
where $v$ is the right null vector of $f_x$, such that $f_xv=0$.

\subsection{First-order LDT Approximation}\label{sec:1st_order}
The first-order approximation replaces the nonlinear bifurcation boundary~$\mathcal{B}$ with its tangent hyperplane at the instanton~$\lambda^*$, yielding the half-space
\[
\mathcal{H} := \{\lambda : N^\top(\lambda - \lambda^*) \geq 0 \},
\]
and estimates $P(\lambda \in \mathcal{B}^-)$ by the probability of this half-space.
For Gaussian $\lambda \sim \mathcal{N}(\mu, \Sigma)$, this integral reduces to a univariate Gaussian tail (see Theorem~B.2 in \cite{Tong_2021}):
\begin{equation}\label{eq:1st_order_ldt}
P_{\mathrm{1st}}(\lambda \in \mathcal{B}^-) := P(\lambda \in \mathcal{H}) = \gausstail{\sqrt{2I(\lambda^*)}},
\end{equation}
where $\Phi(\cdot)$ denotes the standard normal CDF and $I(\lambda^*)$ is the rate function~\eqref{eq:rate_function_mvnormal} evaluated at the instanton.
Intuitively, $\sqrt{2I(\lambda^{*})}$ is the Mahalanobis distance from the mean~$\mu$ to the boundary, and $\Phi(-\cdot)$ converts that distance into a tail probability.
This is identical to the First-Order Reliability Method (FORM) with reliability index $\beta = \sqrt{2I(\lambda^*)}$; the LDT formalism, however, extends naturally to non-Gaussian distributions such as Gaussian mixtures (Section~\ref{sec:gmm_extension}).

\subsection{Second-order LDT Approximation}\label{sec:2nd_order}
The half-space $\mathcal{H}$ used in the first-order approximation can be a poor fit when the bifurcation boundary is strongly curved near the instanton. The second-order approximation refines this by replacing $\mathcal{B}$ with a local quadratic model that incorporates the curvature of the boundary:
\begin{equation}\label{eq:DG_def}
\mathcal{D}_G := \{\lambda : N^\top (\lambda-\lambda^*)+\tfrac{1}{2}(\lambda-\lambda^*)^\top \II (\lambda-\lambda^*)\ge 0\},
\end{equation}
where $N$ is the normal vector and $\II$ is the second fundamental form of~$\mathcal{B}$ at~$\lambda^*$ (Section~\ref{sec:normal_vector}).
The resulting probability $P(\lambda \in \mathcal{D}_G)$ adjusts the first-order estimate whenever the boundary deviates appreciably from its tangent hyperplane.

Note that, unlike the classical setting in \cite{Tong_2021} where the approximated region is obtained by Taylor-expanding an explicit scalar function~$F$, no such scalar function is available in bifurcation analysis. Instead, we define the approximated region directly through the normal vector~$N$ and the second fundamental form~$\II$ as described in Section~\ref{sec:normal_vector}.

For Gaussian $\lambda\sim\mathcal{N}(\mu,\Sigma)$, we derive a closed-form expression for $P(\lambda\in\mathcal{D}_G)$ following~\cite{Tong_2021}. First, we apply an affine transformation that maps $\lambda$ to a standard normal variable~$\xi$ and simultaneously aligns the instanton direction with the first standard basis vector~$e_1$.
For this standardized configuration, \cite{Tong_2021} derives the complete second-order probability formula, and the general Gaussian case reduces to this canonical form.
Specifically, in the standard-normal coordinate system, the approximated region becomes
\begin{equation}\label{eq:area_D}
\mathcal{D} := \{\xi: e_1^\top (\xi-\xi^*)+\frac{1}{2}(\xi-\xi^*)^\top S(\xi-\xi^*) \ge 0\} \,,
\end{equation}
where $\xi^*=\|\xi^*\|e_1$ is the unique minimizer of $\|\xi\|^2$ on $\mathcal{D}$ and $S$ is the matrix that describes second-order information. 
Following~\cite{Tong_2021}, the approximated probability of $\mathcal{D}$ when $\xi$ follows a standard normal distribution is 
\begin{equation}\label{eq:sn_approx}
\mu_{SN}(\mathcal{D}) \approx \gausstail{\|\xi^*\|}\prod^{n-1}_{i=1}[1-\|\xi^*\|\nu_i(P_nSP_n^\top)]^{-\frac{1}{2}}
\end{equation}
where $P_nSP_n^\top$ denotes the $(n-1)\times (n-1)$ submatrix of $S$ obtained by removing its first row and first column, and $\nu_i(\cdot)$ represents the i-th eigenvalue of the matrix argument.

To apply the standard-normal formula~\eqref{eq:sn_approx} to $P(\lambda \in \mathcal{D}_G)$ with $\lambda \sim \mathcal{N}(\mu,\Sigma)$, we transform $\lambda$ to $\xi$ via the affine mapping,
\begin{equation}\label{eq:affine_transformation}
\lambda = A\xi + \mu, \quad A = \Sigma^{\frac{1}{2}}R,
\end{equation}
where $R$ is a rotation matrix such that $R^\top\Sigma^{-\frac{1}{2}}(\lambda^*-\mu)$ aligns with the first standard basis vector $e_1$, and $\lambda^*$ is the instanton under $\mathcal{N}(\mu, \Sigma)$ on the region $\mathcal{D}_G$.

By this transformation, $\xi$ follows a standard normal distribution. It remains to show that the transformed region $\tilde{\mathcal{D}}_G$ satisfies the requirements of~\eqref{eq:area_D}.

Applying the affine transformation to $\mathcal{D}_G$ yields
\begin{align*}
\tilde{\mathcal{D}}_G
={}&\bigl\{\xi:N^\top(A\xi-A\xi^*)\\
&\quad+\tfrac{1}{2}(A\xi-A\xi^*)^\T\II(A\xi-A\xi^*)\ge0\bigr\}\\
={}&\bigl\{\xi:e_1^\T(\xi-\xi^*)\\
&\quad+\tfrac{1}{2}(\xi-\xi^*)^\T\frac{A^\T\II A}{\|A^\T N\|}(\xi-\xi^*) \ge 0\bigr\}
\end{align*}
Formally, we must verify that $A^\top N$ aligns with $e_1$ for the second equality to hold. This can be done by performing a change of variables from $\lambda$ to $\xi$ in the optimality condition~\eqref{eq:general_optimality}:
\begin{align*}
k^*\nabla_{\xi}\tilde{I}(\xi)A^{-1}-w^{*\top} f_{\lambda}(x,\lambda) = 0,
\end{align*}
and since $A^{-1}$ is invertible and $\tilde{I}(\xi)=\frac{1}{2}\xi^T\xi$, we have 
$$
k^{*}\xi^*=w^{*\top} f_{\lambda}(x,\lambda)A^\top=N A^\top.
$$
Following the change of variable in~\eqref{eq:affine_transformation}, the direction of $\xi^*$ aligns with the first standard unit basis vector $e_1$. Since $k^*\neq 0$, we conclude that $NA^\T$ is in the same direction as $e_1$, and $\tilde{\mathcal{D}}_G$ satisfies the requirements for $\mathcal{D}$.

Finally, we substitute $\|\xi^*\|=\sqrt{2I(\lambda^*)}$ and $S =\frac{A^T\II A}{\|A^TN\|}$ in the standard Gaussian case~\eqref{eq:sn_approx}, to obtain 
\begin{equation}\label{eq:2nd_order_ldt}
\begin{aligned}
P_{\mathrm{2nd}}(\lambda\in \mathcal{B}^-)
&= P_{\mathrm{1st}}(\lambda \in \mathcal{B}^-)
\\
&\quad\times
\prod^{n-1}_{i=1}
\left[
1-\frac{\sqrt{2I(\lambda^*)}}{\|A^TN\|}\nu_i(P_nA^T\II AP_n^T)
\right]^{-\frac{1}{2}}
\end{aligned}
\end{equation}
where $\II$ is defined in Section~\ref{sec:normal_vector}. This yields the second-order LDT approximation for the voltage collapse region.

\begin{figure}[t]
  \centering
  \includegraphics[width=\linewidth]{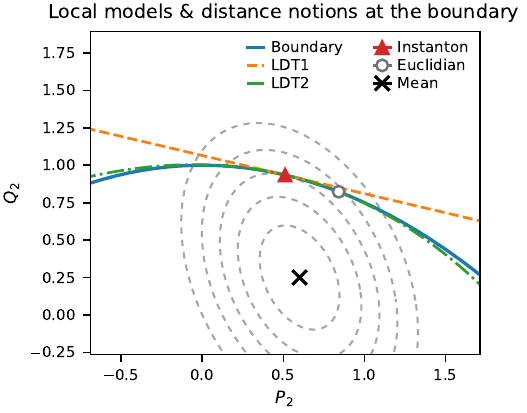}
  \caption{%
  Local geometry at the bifurcation boundary for the 2-bus system.
  True boundary (blue), LDT-1 tangent (orange, dashed), and LDT-2 quadratic (green, dash-dot)
  at the instanton $\lambda^*$ (red triangle). Gray dashed ellipses are level sets of
  $I(\lambda)$ for $N(\mu,\Sigma)$. The separation highlights how anisotropy/rotation
  in $\Sigma$ changes the most probable failure direction and why LDT-2 improves LDT-1 by
  accounting for boundary curvature.
  }
  \label{fig:local-ldt-geometry}
\end{figure}

The derivations above admit a concise geometric interpretation.
For Gaussian uncertainty, the rate function is the Mahalanobis distance
\eqref{eq:rate_function_mvnormal}.
The instanton $\lambda^*$ obtained from the KKT system is precisely the
\emph{Mahalanobis-closest} point on the bifurcation surface; equivalently,
it is the point where a Mahalanobis level set of $I$ is tangent to the boundary.
The first-order LDT approximation replaces the nonlinear boundary by its
\emph{tangent hyperplane} at $\lambda^*$, whereas the second-order approximation uses the
\emph{local quadratic} model determined by the second fundamental form (curvature)
to correct the prefactor.
By contrast, the closest bifurcation point minimizes the
Euclidean distance from $\mu$ to the boundary and generally does
not coincide with $\lambda^*$ when $\Sigma$ is anisotropic or rotated.
Figure~\ref{fig:local-ldt-geometry} illustrates these ideas.
\vspace{2mm}

\subsection{Gaussian Mixture Extension}\label{sec:gmm_extension}
The preceding subsections derived probability estimates under a single Gaussian distribution.
We now extend the framework to Gaussian mixture models (GMMs), which can capture multimodal or skewed uncertainty patterns that a single Gaussian cannot represent.
Suppose that
\begin{equation}
\lambda \sim \sum_{i=1}^{M} \pi_i \mathcal{N}(\mu_i,\Sigma_i), \quad
\pi_i > 0, \quad \sum_{i=1}^{M}\pi_i = 1.
\end{equation}
Following \cite{tong2022optimization}, the cumulant generating function of $\lambda$ is
\begin{equation}
S(\eta)=\log\left(\sum_{i=1}^{M}\pi_i
\exp\left(\eta^\top\mu_i+\frac{1}{2}\eta^\top\Sigma_i\eta\right)\right),
\end{equation}
and the rate function $I = S^*$ (the convex conjugate) no longer admits a closed-form expression.
Instead, the instanton is computed by introducing the dual variable $\eta=\nabla I(\lambda)$ and working with $S$ directly.
Introducing $\eta$ as an additional unknown, the instanton system~\eqref{eq:kkt_system} reduces to
\begin{subequations}\label{eq:gmm_kkt_system}
\begin{align}
f(x,\lambda) &= 0 \,, \label{eq:gmm_power_flow}\\
w^\top f_x(x,\lambda) &= 0 \,, \label{eq:gmm_bifurcation}\\
k \left(w^\top f_\lambda(x,\lambda)\right)^\top &= \eta \,, \\
\frac{\sum_{i=1}^{M}\rho_i(\eta)\left(\mu_i+\Sigma_i\eta\right)}
{\sum_{i=1}^{M}\rho_i(\eta)}
&= \lambda \,, \label{eq:gmm_gradient_relation}\\
w^\top f_\lambda(x,\lambda)\left(w^\top f_\lambda(x,\lambda)\right)^\top &= 1 \label{eq:gmm_normalization}\,,
\end{align}
\end{subequations}
where $\rho_i(\eta)=\pi_i\exp\left(\eta^\top\mu_i+\frac{1}{2}\eta^\top\Sigma_i\eta\right)$.
Equation~\eqref{eq:gmm_gradient_relation} is the gradient relation $\lambda = \nabla S(\eta)$, which maps the dual variable back to the primal space; it replaces the explicit gradient~\eqref{eq:rate_gradient} used in the Gaussian case.
The power-flow and bifurcation conditions~\eqref{eq:gmm_power_flow}--\eqref{eq:gmm_bifurcation} and the normalization~\eqref{eq:gmm_normalization} are unchanged from~\eqref{eq:kkt_system}.
Crucially, the geometric quantities $N=w^\top f_\lambda$ and $\II$ depend only on the bifurcation surface and are computed exactly as before.

Because a Gaussian mixture is a weighted sum of its components, the probability of any measurable set decomposes as $P(\lambda \in \mathcal{A}) = \sum_i \pi_i P_i(\lambda \in \mathcal{A})$, so each LDT estimate likewise decomposes componentwise.
Using the formulas in \cite{tong2022optimization}, the first-order approximation becomes
\begin{equation}\label{eq:gmm_1st_order_ldt}
P_{\mathrm{1st}}^{\mathrm{GM}}(\lambda \in \mathcal{B}^-)
= \sum_{i=1}^{M}\pi_i
\Phi\left(
-\frac{N^\top(\lambda^*-\mu_i)}
{\|\Sigma_i^{\frac{1}{2}}N\|}
\right).
\end{equation}
For the second-order approximation, each component requires its own tangency point on the quadratic model~\eqref{eq:DG_def}:
\begin{equation*}
\tilde{\lambda}_i := \arg\min_{\lambda \in \partial \mathcal{D}_G}
\tfrac{1}{2}\|\lambda-\mu_i\|_{\Sigma_i^{-1}}^2,
\end{equation*}
that is, $\tilde{\lambda}_i$ is the Mahalanobis-closest point on the quadratic boundary to the $i$-th component mean~$\mu_i$.
The local outward normal at $\tilde{\lambda}_i$ is
$
N_i = N + \II(\tilde{\lambda}_i-\lambda^*).
$
This gives
\begin{equation}\label{eq:gmm_2nd_order_ldt}
P_{\mathrm{2nd}}^{\mathrm{GM}}(\lambda \in \mathcal{B}^-)
= \sum_{i=1}^{M}\pi_i P_{2,i},
\end{equation}
where $P_{2,i}$ is exactly the Gaussian second-order formula \eqref{eq:2nd_order_ldt} evaluated with $(\mu,\Sigma,\lambda^*,N)$ replaced by $(\mu_i,\Sigma_i,\tilde{\lambda}_i,N_i)$. Hence the same normal-vector and curvature calculations can be reused componentwise for each Gaussian component.

\subsection{Implementation for Large-scale Systems}

To solve \eqref{eq:kkt_system}, we apply Newton's method to the full residual map
\[
G(z)=0, \qquad z=(x,\lambda,w,k),
\]
so that each iteration has the form
\[
DG(z^{(\ell)})\,\Delta z^{(\ell)}=-G(z^{(\ell)}), 
\qquad 
z^{(\ell+1)}=z^{(\ell)}+\Delta z^{(\ell)}.
\]
Thus, the dominant cost in the instanton computation is the repeated solution of large sparse linear systems involving the derivative blocks induced by \(f_x\) and \(f_\lambda\).

Once the instanton \((x^*,\lambda^*,w^*,k^*)\) has been computed, the remaining task is the evaluation of the geometric quantities entering \eqref{eq:2nd_order_ldt}. The key intermediate quantity is \(x_\lambda\). Since \(f_x(x^*,\lambda^*)\) is singular at the bifurcation point, \(x_\lambda\) is computed columnwise by solving the bordered sparse systems
\begin{equation}
\mathcal{M}^* :=
\begin{bmatrix}
f_x(x^*,\lambda^*) & v^* \\
w^{*\top} f_{xx}(x^*,\lambda^*)(\,\cdot\,,v^*) & 0
\end{bmatrix}.
\end{equation}
Then, for each \(j=1,\dots,m\),
\begin{equation}
\mathcal{M}^*
\begin{bmatrix}
x_{\lambda_j}\\
\alpha_j
\end{bmatrix}
=
\begin{bmatrix}
-\,f_{\lambda_j}(x^*,\lambda^*)\\
0
\end{bmatrix}.
\end{equation}
where \(v^*\) is the right null vector of \(f_x(x^*,\lambda^*)\), \(f_{\lambda_j}\) denotes the \(j\)-th column of \(f_\lambda\), and \(\alpha_j\) is an auxiliary scalar. This avoids any dense inversion at the singular point and yields the derivatives required to assemble the quantities in Section~\ref{sec:normal_vector}. In the power-flow setting considered here, \(f_{x\lambda}=0\) and \(f_{\lambda\lambda}=0\), so the only second derivatives needed are those contained in \(f_{xx}\).

%\section{Computational implementation}

\section{Experiments}

To demonstrate the methodology, we perform experiments on two-bus and five-bus systems described in \cite{dobsondirect1991}.

\subsection{Two-bus system}
With the two-bus system as an example, we show how the computation of the LDT probabilities is affected by different uncertainty models, including isotropic and anisotropic Gaussian covariances as well as a Gaussian mixture. The two-bus system can be described with the following equations
\begin{align*}
-4V\sin\alpha-\lambda_1&=0 \, \\
-4V^2+4V\cos\alpha-\lambda_2&=0 \,.
\end{align*}
where $V$ is the voltage of the load bus and $\alpha$ its angle. The load injections are modeled by the random vector $(\lambda_1, \lambda_2)$. One of the benefits of considering this simple system is that the bifurcation surface can be derived analytically \cite{dobsondirect1991}
\begin{equation}
    \lambda_1^2+4\lambda_2-4=0 \,.
\end{equation}

\subsubsection{Gaussian distribution, isotropic covariance}
When $\lambda$ is normally distributed with identity covariance, computing the instanton is equivalent to finding the closest bifurcation point to $\mu$ in the voltage-collapse region, since the rate function reduces to the squared Euclidean distance from $\mu$.

Indeed, we observe that when computing the instanton using \eqref{eq:kkt_system} and using mean vector (nominal point) $\mu=[0.5,0.3]$, the instanton is found at $[0.703,0.877]$ which aligns with the closest bifurcation point obtained in \cite{dobson1993}.

To visualize the shape of the bifurcation surface in the $\lambda$-space, we generated $10^{5}$ random samples of $\lambda$. Since the 2-bus bifurcation boundary is available analytically as
\[
\lambda_1^2 + 4\lambda_2 - 4 = 0,
\]
we classify a sample as collapsed whenever
\[
\lambda_1^2 + 4\lambda_2 - 4 > 0,
\]
and as solvable otherwise. 

In Figure~\ref{fig:local-ldt-geometry},  we can observe that the bifurcation boundary is curving outward from the mean. This means that the second-order LDT approximation should yield a larger probability than the first-order approximation, since the curvature correction factor will exceed one. We verify this fact in the following experiments.

\subsubsection{Gaussian distribution with general covariance matrix}
In this experiment we consider a more general covariance matrix, $\Sigma_{0}=\mathrm{diag}\{0.6,1.0\}$.
To test the performance of the approximation as the probability of voltage collapse decreases, we scale the covariance by a logarithmically spaced sequence $c_i$ ranging from $0.631$ to $0.0158$ and consider the Gaussian family $\mathcal{N}(\mu,c_i\Sigma_0)$ with mean vector $\mu = [0.5,0.3]$. 

For each $c_i$, the instanton is found by solving \eqref{eq:kkt_system} after substituting \eqref{eq:rate_gradient}. The resulting tangency point is then used in the first-order and second-order LDT formulas \eqref{eq:1st_order_ldt} and \eqref{eq:2nd_order_ldt}. Reference probabilities are computed by direct Monte Carlo when $c_i \geq 1.227\times 10^{-1}$ and by Gaussian importance sampling for smaller covariance scales, using a proposal centered at the instanton. In the reference calculation we use $4\times 10^5$ direct Monte Carlo samples and $1.5\times 10^5$ importance-sampling samples; these counts keep the estimator variance small while remaining inexpensive in the 2-bus setting.

\begin{figure}[t]
    \centering
    \includegraphics[width=\linewidth]{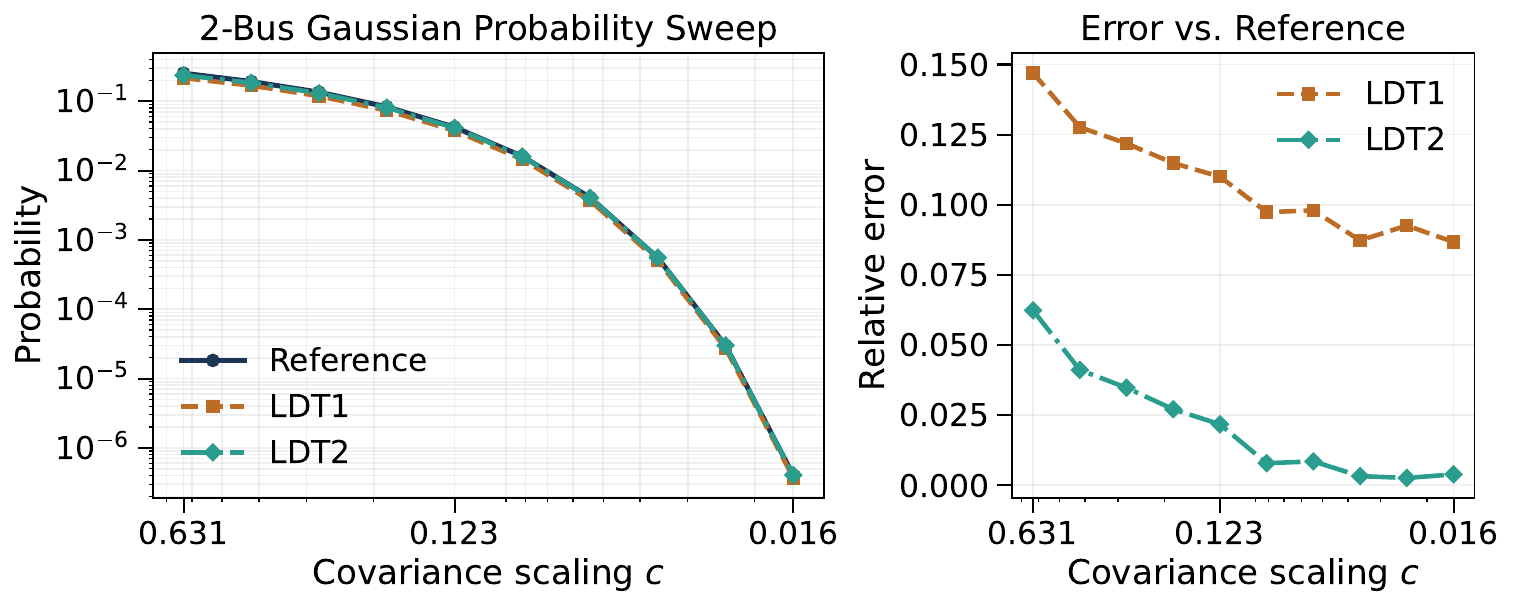}
    \caption{Two-bus probability sweep for Gaussian uncertainty. The reference curve is computed by direct Monte Carlo for the larger covariance scales and by importance sampling in the rare-event regime.}
    \label{fig:gaussian_2bus_sweep}
\end{figure}

\begin{table}[t]
\centering
\caption{Comparison of Reference and LDT Approximations under Covariance Scaling in the 2-bus System}
\label{tab:ldt_2bus_results}
\resizebox{\columnwidth}{!}{\input{gaussian_2bus_table.tex}}
\end{table}

Figure~\ref{fig:gaussian_2bus_sweep} and Table~\ref{tab:ldt_2bus_results} show that the exact Gaussian first- and second-order formulas are already accurate across the full sweep, with the second-order correction improving the first-order estimate at every tested covariance scale. Over the full sweep, the maximum relative error decreases from $14.7\%$ for LDT1 to $6.2\%$ for LDT2. Once the reference probability drops below $10^{-3}$, the second-order approximation stays within $0.4\%$ of the Monte Carlo / importance-sampling reference.

\subsubsection{Gaussian mixture with two components}
To illustrate the non-Gaussian extension in the same 2-bus setting, we consider the two-component Gaussian mixture
\[
\lambda \sim 0.75\,\mathcal{N}(\mu_1,c\Sigma_1)+0.25\,\mathcal{N}(\mu_2,c\Sigma_2),
\]
with
\[
\mu_1=\begin{bmatrix}0.45\\0.25\end{bmatrix},
\qquad
\mu_2=\begin{bmatrix}0.82\\0.52\end{bmatrix},
\]
and
\[
\Sigma_1=\begin{bmatrix}0.60&0\\0&1.00\end{bmatrix},
\qquad
\Sigma_2=\begin{bmatrix}0.35&0.08\\0.08&0.55\end{bmatrix}.
\]
We use the same covariance scaling sweep $c_i$ as in the Gaussian experiment. For each $c_i$, we solve the augmented instanton system \eqref{eq:gmm_kkt_system}, compute the geometric quantities $N$ and $\II$, and then evaluate the componentwise approximations \eqref{eq:gmm_1st_order_ldt} and \eqref{eq:gmm_2nd_order_ldt}. Reference probabilities are obtained by direct Monte Carlo when $c_i \geq 1.227\times 10^{-1}$ and by mixture importance sampling for smaller $c_i$, using proposal components centered at the first-order tangency points. For these references we use $10^6$ direct Monte Carlo samples and $4\times 10^5$ importance-sampling samples.

\begin{figure}[t]
    \centering
    \includegraphics[width=\linewidth]{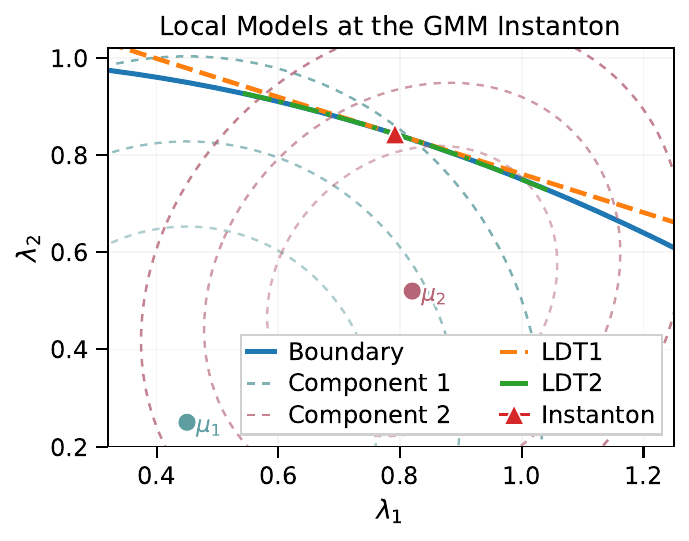}
    \caption{Local geometry at the GMM instanton in the 2-bus system.}
    \label{fig:gmm_2bus_geometry}
\end{figure}

Figure~\ref{fig:gmm_2bus_geometry} plays the same role as Fig.~\ref{fig:local-ldt-geometry} in the Gaussian case: it shows that the local boundary geometry is unchanged, but the uncertainty model is now multimodal, so the probability evaluation must be carried out componentwise.

\begin{figure}[t]
    \centering
    \includegraphics[width=\linewidth]{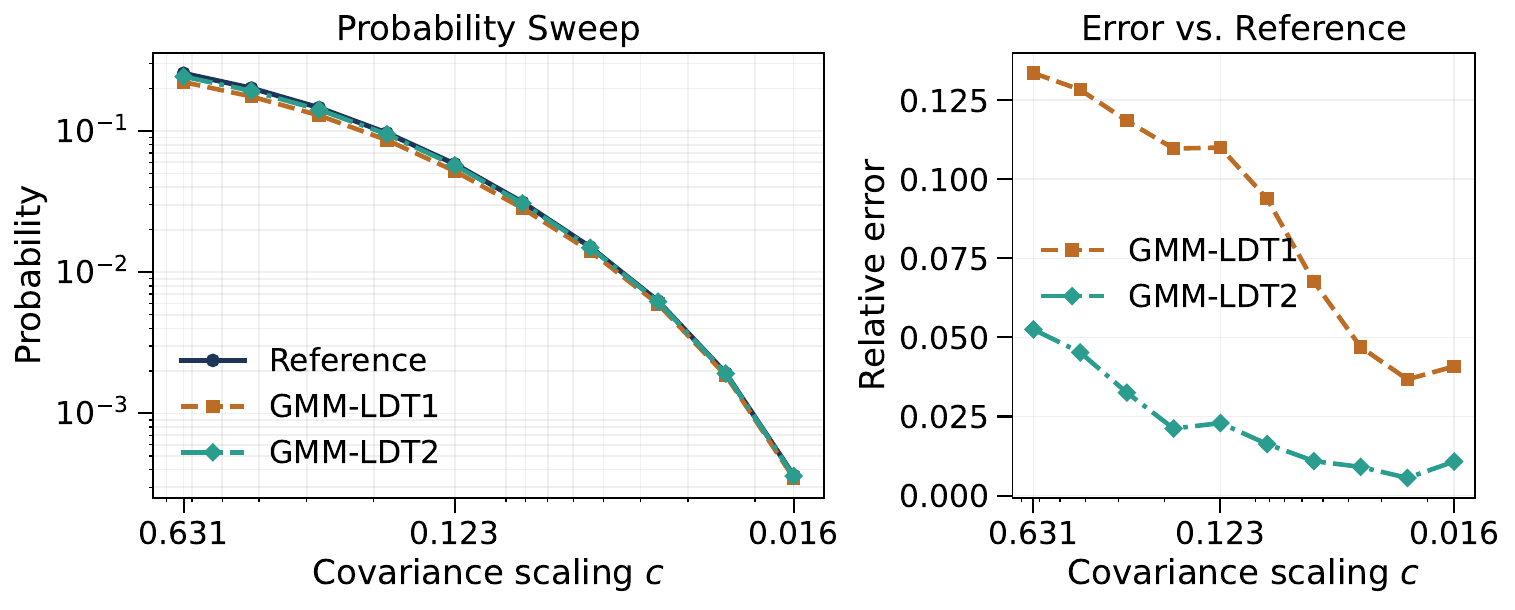}
    \caption{Two-bus probability sweep for the Gaussian-mixture extension. The reference curve is computed by Monte Carlo or importance sampling, while GMM-LDT1 and GMM-LDT2 correspond to \eqref{eq:gmm_1st_order_ldt} and \eqref{eq:gmm_2nd_order_ldt}.}
    \label{fig:gmm_2bus_sweep}
\end{figure}

\begin{table}[t]
\centering
\caption{Comparison of Reference and GMM-LDT Approximations under Covariance Scaling in the 2-bus System}
\label{tab:gmm_2bus_results}
\resizebox{\columnwidth}{!}{\input{gmm_2bus_table.tex}}
\end{table}

Table~\ref{tab:gmm_2bus_results} and Fig.~\ref{fig:gmm_2bus_sweep} show that the second-order correction remains beneficial in the Gaussian-mixture setting. Over the full sweep, the maximum relative error drops from $13.4\%$ for GMM-LDT1 to $5.3\%$ for GMM-LDT2. In the rarer regime where the reference probability is below $10^{-2}$, the second-order approximation stays within $1.1\%$ of the reference values, confirming that the Gaussian geometry machinery can be reused componentwise for a multimodal uncertainty model.

\subsection{Five-bus system}
We now apply the computation to a 5-bus network to demonstrate its scalability and accuracy in higher dimensions. The system is defined by six uncertain load parameters $\lambda = (P_2, Q_2, P_4, Q_4, P_5, Q_5)$ and admittance matrix described in Table~\ref{tab:ybus_5bus}. For the Gaussian experiments, the stochastic loading vector is modeled with mean $\mu = [1.15, 0.6, 0.7, 0.3, 0.7, 0.4]$ and covariance $c_i\Sigma_0$, where $\Sigma_0$ is a generated positive semi-definite matrix. We also consider a Gaussian-mixture uncertainty model. In both cases, reference probabilities are computed by direct Monte Carlo in the moderate-probability regime and by importance sampling for the rarer cases.

\begin{table}[h]
\centering
\caption{Admittance Matrix $Y_{\mathrm{bus}}$ for the 5-bus System}
\label{tab:ybus_5bus}
\renewcommand{\arraystretch}{1.15}
\setlength{\tabcolsep}{5pt}
\resizebox{\columnwidth}{!}{
\begin{tabular}{c|ccccc}
\hline
\textbf{Bus} & \textbf{1} & \textbf{2} & \textbf{3} & \textbf{4} & \textbf{5} \\ 
\hline
\textbf{1} & $3.241 - j13.085$ & $-1.401 - j5.602$ & $0$ & $0$ & $-1.841 - j7.484$ \\
\textbf{2} & $-1.401 - j5.602$ & $3.242 - j12.486$ & $-1.841 - j7.484$ & $0$ & $0$ \\
\textbf{3} & $0$ & $-1.841 - j7.484$ & $3.671 - j14.762$ & $-0.700 - j2.801$ & $-1.130 - j4.477$ \\
\textbf{4} & $0$ & $0$ & $-0.700 - j2.801$ & $1.634 - j6.236$ & $-0.934 - j3.435$ \\
\textbf{5} & $-1.841 - j7.484$ & $0$ & $-1.130 - j4.477$ & $-0.934 - j3.435$ & $3.905 - j15.396$ \\
\hline
\end{tabular}
}
\end{table}

\subsubsection{Gaussian distribution with general covariance matrix}

As we did before, to investigate the asymptotic behavior of the probability of voltage collapse, we scale the covariance matrix by a factor $c_i$ to construct a family of Gaussian distributions 
$\mathcal{N}(\mu,\, c_i \Sigma_0)$, 
with $c_i \in [0.64,\, 0.02]$. 
For each scaling factor $c_i$, we compute the corresponding instanton $\lambda^*$, 
the normal vector $N = w^\top f_\lambda$, 
the second fundamental form $\mathrm{II}$, 
and the associated first-order and second-order LDT probability approximations. 
Due to the increased size of the 5-bus system and the rapid variation of the collapse probability under shrinking covariance, 
we employ direct Monte Carlo when $c_i \geq 6.3\times 10^{-2}$ and Gaussian importance sampling below that threshold, with proposal distributions centered at the corresponding instantons. This gives a stable reference curve across both the moderate and rare-event regimes. In the reference calculation we use $4000$ direct Monte Carlo samples and $3000$ importance-sampling samples; because each sample requires a nonlinear power-flow solve, we chose these counts to balance runtime against estimator variance, with importance sampling preventing the rare-event points from becoming too noisy.

\begin{figure}[t]
    \centering
    \includegraphics[width=\linewidth]{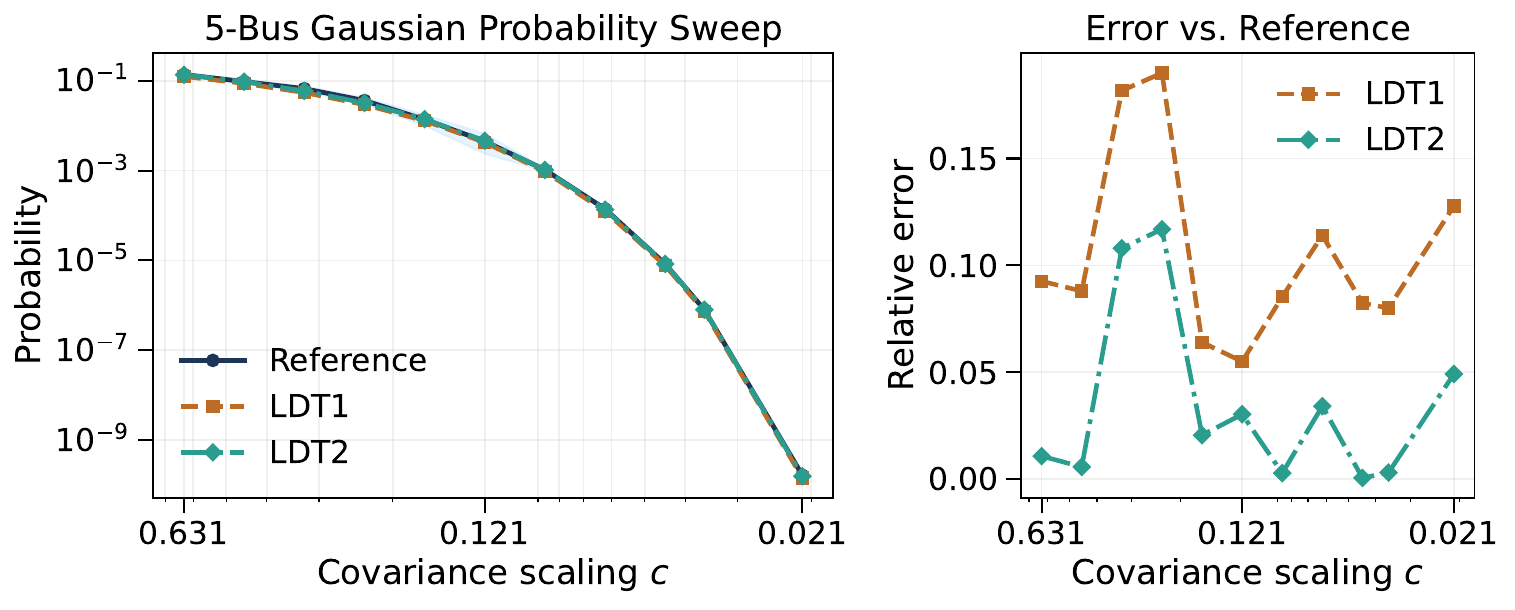}
    \caption{Five-bus probability sweep for Gaussian uncertainty. The reference curve combines direct Monte Carlo with importance sampling, while LDT1 and LDT2 denote the first-order and second-order approximations.}
    \label{fig:gaussian_5bus_sweep}
\end{figure}

\begin{table}[t]
\centering
\caption{Comparison of Reference and LDT Approximations under Covariance Scaling in the 5-bus System}
\label{tab:5bus_covariance_scaling}
\resizebox{\columnwidth}{!}{\input{gaussian_5bus_table.tex}}
\end{table}

Figure~\ref{fig:gaussian_5bus_sweep} and Table~\ref{tab:5bus_covariance_scaling} show the same pattern in higher dimension: the second-order correction is usually closer to the reference and becomes especially useful deep in the rare-event regime. After switching to importance sampling from $c=8.660\times 10^{-2}$ downward, the relative-error curves become much smoother. For reference probabilities below $10^{-3}$ the maximum relative error drops from $12.8\%$ for LDT1 to $4.9\%$ for LDT2, and at $c=8.660\times 10^{-2}$ the second-order estimate is within $0.3\%$ of the importance-sampling reference.

\subsubsection{Gaussian mixture with two components}
To check that the Gaussian-mixture extension also remains practical in higher dimension, we repeated the 5-bus experiment with
\[
\lambda \sim 0.8\,\mathcal{N}(\mu_1,c\Sigma_1)+0.2\,\mathcal{N}(\mu_2,c\Sigma_2),
\]
where
\[
\mu_1=\begin{bmatrix}1.15\\0.60\\0.70\\0.30\\0.70\\0.40\end{bmatrix}, \qquad
\mu_2=\begin{bmatrix}1.33\\0.66\\0.86\\0.35\\0.85\\0.46\end{bmatrix},
\]
$\Sigma_1=\mathrm{diag}(1.0,0.8,0.3,1.0,1.0,1.0)$, and $\Sigma_2$ is a correlated positive definite matrix constructed in the Python verification script. For each $c$, we solve \eqref{eq:gmm_kkt_system}, compute $N$ and $\II$ exactly as in the Gaussian case, and then evaluate the componentwise GMM-LDT1 and GMM-LDT2 approximations. Reference probabilities are obtained by direct Monte Carlo for the two largest covariance scales and by mixture importance sampling for the remaining three rare-event cases. Here we use $6000$ direct Monte Carlo samples and $4000$ mixture importance-sampling samples.

For the 5-bus Gaussian workflow, representative runtimes on the machine used for these experiments are about $9.3$ s for direct Monte Carlo with $4000$ samples, $11.1$ s for importance sampling with $3000$ samples, $0.28$ s for LDT1, and $0.33$ s for LDT2. Thus the dominant cost in the reference calculations is repeated power-flow solution, while the second-order LDT correction adds only a modest overhead beyond the instanton solve already needed for LDT1.

\begin{figure}[t]
    \centering
    \includegraphics[width=\linewidth]{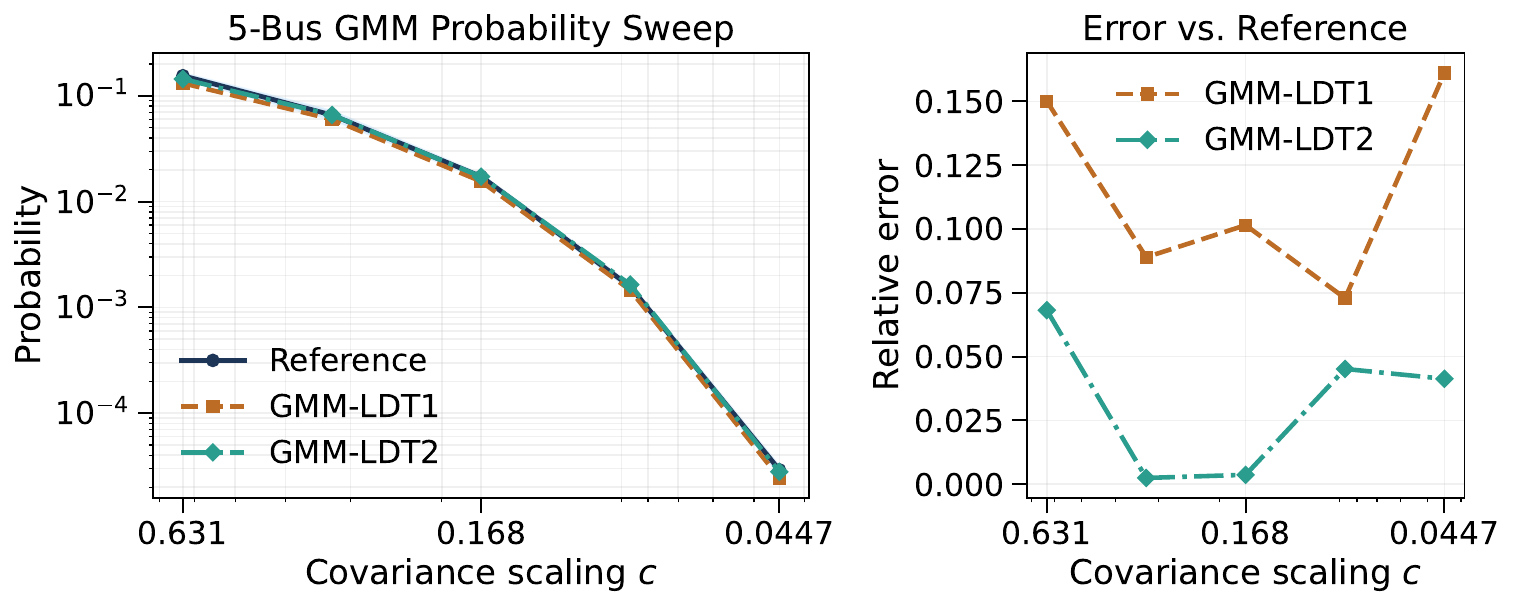}
    \caption{Five-bus probability sweep for the Gaussian-mixture extension. The second-order approximation remains uniformly closer to the Monte Carlo / importance-sampling reference than the first-order approximation over the tested covariance scales.}
    \label{fig:gmm_5bus_sweep}
\end{figure}

\begin{table}[t]
\centering
\caption{Comparison of Reference and GMM-LDT Approximations under Covariance Scaling in the 5-bus System}
\label{tab:gmm_5bus_results}
\resizebox{\columnwidth}{!}{\input{gmm_5bus_table.tex}}
\end{table}

Table~\ref{tab:gmm_5bus_results} and Fig.~\ref{fig:gmm_5bus_sweep} show that the second-order correction remains useful in the multimodal setting. Across the full sweep, the maximum relative error decreases from $16.1\%$ for GMM-LDT1 to $6.8\%$ for GMM-LDT2. In the rarer cases with reference probability at or below $1.7\times 10^{-2}$, the second-order approximation stays within $4.6\%$ of the reference values. This gives a direct 5-bus validation that the Gaussian-mixture extension can reuse the same geometric computations developed for the Gaussian model.

\subsection{Computational Performance}

We evaluated the methodology on the IEEE 14, IEEE 57, and IEEE 118 test cases. In all cases, uncertainty was restricted to loads at five PQ buses, so that the parameter dimension was fixed at \(m=10\). The implementation used the Julia programming language and the PowerModels.jl library \cite{8442948}, and the experiments were performed on a 2024 Apple MacBook Pro M4 Max.

The reported timings consist of three stages: construction of the sparse model, solution of the instanton KKT system, and evaluation of the LDT-2 geometric correction. For each case, one warm-up run was performed in the same Julia session, followed by three measured runs; all reported times are medians over the measured runs.

Table~\ref{tab:workflow_timings} reports the overall workflow timings. Model construction is negligible in all three cases. The dominant costs are the instanton KKT solve and the LDT-2 shape-operator computation. Even for the IEEE 118 case, the total reported cost of the deterministic portion of the workflow is only \(4.26\) seconds.

\begin{table}[h]
\centering
\caption{Workflow timings for the larger-case experiments. All values are medians over three measured runs after one warm-up run. Times are in seconds.}
\label{tab:workflow_timings}
\resizebox{\columnwidth}{!}{
\begin{tabular}{c|cc|ccccc}
\hline
\textbf{Case} & \(\mathbf{n}\) & \(\mathbf{m}\) & \textbf{Build} & \textbf{KKT solve} & \textbf{LDT-2} & \textbf{Total}\\
\hline
IEEE 14  & 22  & 10 & 0.0013 & 0.0224 & 0.0097 & 0.0330 \\
IEEE 57  & 106 & 10 & 0.0036 & 0.6530 & 0.4597 & 1.1333  \\
IEEE 118 & 181 & 10 & 0.0095 & 2.2633 & 1.9926 & 4.2646  \\
\hline
\end{tabular}}
\end{table}

To identify the source of the LDT-2 cost, Table~\ref{tab:ldt2_breakdown} decomposes the implicit shape-operator computation into its principal stages. The cost of forming \(f_\lambda\), constructing the normal vector and tangent basis, forming the Jacobian matrix \(J_x\), and assembling the fold operator is negligible. The dominant contribution is the family of tangent-direction linear solves required by the implicit differentiation step. This behavior is consistent with the formulation in the previous subsection: once the instanton has been computed, the second-order correction is obtained through a sequence of sparse bordered solves, one for each tangent direction. Since \(m=10\) in all cases, the tangent space has dimension \(m-1=9\), and the average solve time per tangent direction increases substantially with network size.

\begin{table}[h]
\centering
\caption{Breakdown of the LDT-2 shape-operator computation. All values are medians over three measured runs after one warm-up run. Times are in seconds.}
\label{tab:ldt2_breakdown}
\resizebox{\columnwidth}{!}{
\begin{tabular}{c|ccccc|c}
\hline
\textbf{Case}  & \(\mathbf{J_x}\) & \textbf{Fold-op} & \textbf{Tangent solves} & \textbf{Avg.\ solve/dir.} \\
\hline
IEEE 14   & \(3.85\times 10^{-5}\) & \(5.00\times 10^{-7}\) & 0.00963 & 0.00107 \\
IEEE 57    & \(1.62\times 10^{-4}\) & \(8.33\times 10^{-7}\) & 0.45948  & 0.05105 \\
IEEE 118   & \(2.53\times 10^{-4}\) & \(1.08\times 10^{-6}\) & 1.99228  & 0.22136 \\
\hline
\end{tabular}}
\end{table}

The computational bottlenecks are the repeated sparse linear solves arising in the instanton KKT solve and in the tangent-direction solves required by the LDT-2 geometry. In particular, the latter accounts for essentially the entire cost of the shape-operator computation. This, however, can be parallelized to alleviate the computational costs.

\subsection{Comparison with Reliability Methods (FORM/SORM)}
For Gaussian uncertainty, the first-order LDT formula~\eqref{eq:1st_order_ldt} coincides with FORM after identifying the reliability index as $\beta=\sqrt{2I(\lambda^\ast)}$, which is the Mahalanobis distance from the mean to the instanton. The second-order formula~\eqref{eq:2nd_order_ldt} has the same product structure as SORM: in the standard-normal coordinate system obtained via the affine transformation~\eqref{eq:affine_transformation}, the eigenvalues $\nu_i(P_nA^\top\II AP_n^\top)/\|A^\top N\|$ play the role of the principal curvatures of the limit-state surface, and the resulting correction factor matches the classical SORM prefactor \cite{Du2001,Tong_2021}. In this sense, the Gaussian results in Sections~\ref{sec:1st_order} and~\ref{sec:2nd_order} can also be read as FORM/SORM results expressed in the LDT language of rate functions and instantons.

The advantage of the present framework is that it extends naturally to non-Gaussian uncertainties through the same instanton formulation. As shown in the Gaussian-mixture extension (Section~\ref{sec:gmm_extension}), the bifurcation-surface quantities, namely the instanton, normal vector, and second fundamental form, are computed exactly as before, while only the probability evaluation changes to a componentwise sum over the mixture components.

\section{Conclusion and Future research}
Our work establishes a rigorous framework for estimating the probability of voltage collapse by leveraging the principles of large deviation theory (LDT). By reformulating the problem of voltage instability as a rare-event probability estimation, we implemented both first-order and second-order LDT approximations that capture the asymptotic behavior of collapse probabilities under random loading. This framework unifies previous approaches, offering a generalization of existing methods such as the Mahalanobis-distance-based formulation of Mittelstaedt et al. \cite{Mittelstaedt2017}. Numerical experiments on 2-bus and 5-bus systems demonstrate that the LDT-based approach achieves acceptable accuracy in the rare-event regime where standard Monte Carlo methods become computationally infeasible, since crude Monte Carlo requires on the order of $1/p$ samples to resolve an event of probability $p$ with fixed relative accuracy.

The next step is to scale the method to networks with hundreds or thousands of buses. In practice, we can combine stable iterative solvers \cite{dobson1993, Mittelstaedt2017} with our LDT framework, and use continuation from a known base point \cite{ContinuationPF} together with trust regions or Levenberg-Marquardt damping \cite{levenberg2022} to avoid divergence near the bifurcation boundary. These ingredients provide warm starts and local geometric information for the instanton search, making it more reliable and faster on large systems.

While most of the experiments in this paper focus on Gaussian uncertainty, the Gaussian-mixture derivation together with the 2-bus and 5-bus validations above show that the framework is not limited to that case. We know that load uncertainty can be skewed or heavy-tailed as discussed in \cite{NonGaussianity2018}. Extending the same verification pipeline to larger networks under richer non-Gaussian models is a natural next step.

\bibliographystyle{IEEEtran}
\bibliography{biblio}

\end{document}

%% file: gaussian_2bus_table.tex
\begin{tabular}{c|c|ccc}
\hline
\textbf{Scaling factor $c$} & \textbf{Ref.} & \textbf{Reference} & \textbf{LDT1} & \textbf{LDT2} \\
\hline
$6.310e-01$ & MC & $2.536e-01$ & $2.163e-01$ & $2.378e-01$ \\
$4.190e-01$ & MC & $1.923e-01$ & $1.678e-01$ & $1.844e-01$ \\
$2.783e-01$ & MC & $1.351e-01$ & $1.187e-01$ & $1.304e-01$ \\
$1.848e-01$ & MC & $8.306e-02$ & $7.352e-02$ & $8.081e-02$ \\
$1.227e-01$ & MC & $4.222e-02$ & $3.757e-02$ & $4.130e-02$ \\
$8.149e-02$ & IS & $1.606e-02$ & $1.449e-02$ & $1.593e-02$ \\
$5.412e-02$ & IS & $4.086e-03$ & $3.686e-03$ & $4.052e-03$ \\
$3.594e-02$ & IS & $5.524e-04$ & $5.042e-04$ & $5.542e-04$ \\
$2.387e-02$ & IS & $3.012e-05$ & $2.733e-05$ & $3.004e-05$ \\
$1.585e-02$ & IS & $4.034e-07$ & $3.684e-07$ & $4.050e-07$ \\
\hline
\end{tabular}

%% file: gmm_2bus_table.tex
\begin{tabular}{c|c|ccc}
\hline
\textbf{Scaling factor $c$} & \textbf{Ref.} & \textbf{Reference} & \textbf{GMM-LDT1} & \textbf{GMM-LDT2} \\
\hline
$6.310e-01$ & MC & $2.563e-01$ & $2.221e-01$ & $2.429e-01$ \\
$4.190e-01$ & MC & $2.014e-01$ & $1.755e-01$ & $1.923e-01$ \\
$2.783e-01$ & MC & $1.462e-01$ & $1.289e-01$ & $1.415e-01$ \\
$1.848e-01$ & MC & $9.690e-02$ & $8.627e-02$ & $9.485e-02$ \\
$1.227e-01$ & MC & $5.842e-02$ & $5.200e-02$ & $5.708e-02$ \\
$8.149e-02$ & IS & $3.129e-02$ & $2.836e-02$ & $3.078e-02$ \\
$5.412e-02$ & IS & $1.507e-02$ & $1.405e-02$ & $1.491e-02$ \\
$3.594e-02$ & IS & $6.245e-03$ & $5.952e-03$ & $6.189e-03$ \\
$2.387e-02$ & IS & $1.928e-03$ & $1.857e-03$ & $1.917e-03$ \\
$1.585e-02$ & IS & $3.640e-04$ & $3.491e-04$ & $3.601e-04$ \\
\hline
\end{tabular}

%% file: gaussian_5bus_table.tex
\begin{tabular}{c|c|ccc}
\hline
\textbf{Scaling factor $c$} & \textbf{Ref.} & \textbf{Reference} & \textbf{LDT1} & \textbf{LDT2} \\
\hline
$6.310e-01$ & MC & $1.378e-01$ & $1.250e-01$ & $1.363e-01$ \\
$4.532e-01$ & MC & $9.575e-02$ & $8.733e-02$ & $9.521e-02$ \\
$3.255e-01$ & MC & $6.675e-02$ & $5.461e-02$ & $5.954e-02$ \\
$2.337e-01$ & MC & $3.625e-02$ & $2.936e-02$ & $3.201e-02$ \\
$1.679e-01$ & MC & $1.375e-02$ & $1.287e-02$ & $1.403e-02$ \\
$1.206e-01$ & MC & $4.500e-03$ & $4.253e-03$ & $4.636e-03$ \\
$8.660e-02$ & IS & $1.039e-03$ & $9.508e-04$ & $1.037e-03$ \\
$6.219e-02$ & IS & $1.400e-04$ & $1.240e-04$ & $1.352e-04$ \\
$4.467e-02$ & IS & $8.366e-06$ & $7.677e-06$ & $8.370e-06$ \\
$3.600e-02$ & IS & $7.955e-07$ & $7.319e-07$ & $7.979e-07$ \\
$2.100e-02$ & IS & $1.644e-10$ & $1.434e-10$ & $1.563e-10$ \\
\hline
\end{tabular}

%% file: gmm_5bus_table.tex
\begin{tabular}{c|c|ccc}
\hline
\textbf{Scaling factor $c$} & \textbf{Ref.} & \textbf{Reference} & \textbf{GMM-LDT1} & \textbf{GMM-LDT2} \\
\hline
$6.310e-01$ & MC & $1.548e-01$ & $1.316e-01$ & $1.443e-01$ \\
$3.255e-01$ & MC & $6.583e-02$ & $5.998e-02$ & $6.599e-02$ \\
$1.679e-01$ & IS & $1.727e-02$ & $1.552e-02$ & $1.721e-02$ \\
$8.660e-02$ & IS & $1.573e-03$ & $1.458e-03$ & $1.644e-03$ \\
$4.467e-02$ & IS & $2.912e-05$ & $2.443e-05$ & $2.791e-05$ \\
\hline
\end{tabular}